\documentclass[11pt]{amsart}
\usepackage[noBBpl]{mathpazo}
\usepackage{amsfonts}
\usepackage{amsmath}
\usepackage{amssymb}
\usepackage{amsthm}

\input xy
\xyoption {all}

\newtheorem{theorem}{Theorem}
\newtheorem {lemma}{Lemma}

\newtheorem {corollary}{Corollary}
\newtheorem {proposition}{Proposition}

\newcommand \bfw {{\mathbf W}}
\newcommand \bfv {{\mathbf v}}
\newcommand \bfe {{\mathbf E}}
\newcommand \bfs {{\mathbf S}}
\newcommand \bfh {{\mathbf H}}
\newcommand \bfm {{\mathsf m}}

\newcommand \bfg {{\mathsf G}}
\newcommand \bfM {{\mathcal M}}

\newcommand \bfz {{\mathsf Z}}
\newcommand \bfn {{\it n}}

\theoremstyle{definition}
\newtheorem{remark}{Remark}

\input xy
\xyoption {all}

\begin{document}
\title {The Verlinde bundles and the semihomogeneous Wirtinger duality}
\author {Dragos Oprea}
\address {Department of Mathematics}
\address {University of California, San Diego}
\email {doprea@math.ucsd.edu}
\begin{abstract}
We determine the splitting type of the Verlinde vector bundles in higher genus in terms of simple semihomogeneous factors. In agreement with strange duality, the simple factors are interchanged by the Fourier-Mukai transform, and their spaces of sections are naturally dual.  
\end {abstract}

\maketitle
\section{Introduction}

\subsection {The Verlinde bundles} Let $X$ be a smooth complex projective curve of genus $g\geq 1$. The Jacobian $$A=\text {Jac}(X)$$ carries a class of interesting vector bundles, the Verlinde bundles, introduced and studied by Popa \cite {Po}. The fibers of these bundles are the spaces of generalized theta functions of fixed rank and fixed determinant, as the determinant varies in the Jacobian $A$.  

More explicitly, let $U_{X}(r,0)$ denote the moduli space of semistable bundles on $X$ of rank $r$ and degree $0$. There is a morphism $$\det:U_{X}(r,0)\to A$$ taking rank $r$ vector bundles to their determinants. Fixing a theta characteristic $\kappa$ on $X$, we obtain the generalized Theta line bundle $\Theta_{\kappa}$ on $U_{X}(r,0)$ corresponding to the divisorial locus $$\Theta_{\kappa}=\{E: h^{0}(E\otimes \kappa)\neq 0\}\hookrightarrow U_{X}(r,0).$$ 
In particular, when $r=1$, we obtain a symmetric polarization $\Theta_{\kappa}$ on the Jacobian $A$.

The Verlinde bundles are the pushforwards of the pluri-theta bundles on $U_X(r, 0)$ to the Jacobian via the determinant map $$\bfe_{r,k}={\det}_{\star}(\Theta^{k}_{\kappa}).$$ The integer $k$ is called the level, while $r$ will somewhat misleadingly be referred to as the rank.  Henceforth, the subscript $\kappa$ decorating the theta bundles will be suppressed for simplicity of notation. The same reference theta characteristic will be used throughout. 

The goal of this paper is to express the Verlinde bundles $\bfe_{r,k}$ in terms of indecomposable factors which we proceed to describe. 

\subsection {The semihomogeneous vector bundles} More generally, we will consider a class of vector bundles over abelian varieties constructed by Atiyah in dimension $1$ \cite {A}, and by Mukai in all dimensions \cite {mukai}. To set the stage, let $(A, \Theta)$ be a principally polarized complex abelian variety with $\Theta$ symmetric. For any {\it odd} coprime positive integers $(a, b)$, we will show that there is a unique simple symmetric semihomogeneous bundle $\bfw_{a,b}$ on $A$ such that $$\text{rank }\bfw_{a,b}=a^{g} \text { and } \det \bfw_{a,b}=\Theta^{a^{g-1}b}.$$ (The bundles $\bfw_{a,b}$ may also be defined for even values of $a$ or $b$ by equation \eqref{heisenberg} below.) Semihomogeneity is the requirement that for all $x\in A$, there exists a line bundle $y$ on $A$ such that $$t_{x}^{\star}\bfw_{a,b}=\bfw_{a,b}\otimes y.$$ Here, $t_{x}:A\to A$ denotes the translation by $x$ on the abelian variety. Symmetry means that $$(-1)^{\star}\bfw_{a,b}\cong \bfw_{a,b}.$$ 

\subsection {The splitting of the Verlinde bundles} When $\gcd(r,k)=1$, it is not difficult to show that \begin{equation}\label{coprime1}\bfe_{r, k}=\bigoplus\bfw_{r,k}.\end{equation} The splitting of the Verlinde bundles in arbitrary ranks and levels is less straightforward. An explicit decomposition will be presented below. 

\subsubsection{} Define $$\bfv_{g}({r,k})=(r+k)^{r(g-1)}\sum_{\stackrel{S\subset\{1, \ldots, r+k\}}{ |S|=r}}\prod_{s\neq t \in S}\left|2\sin \frac {s-t}{r+k}\pi\right|^{1-g}.$$ It is easy to observe the level-rank symmetry \begin{equation}\label{symmetry1}\bfv_{g}(r,k)=\bfv_g(k,r)\end{equation} obtained by replacing the set $S$ in the above sum by its complement. 

The numbers $\bfv_g(r,k)$ are related to the Verlinde dimensions: $$h^{0}(SU_{X}(r), \mathcal L^{k})=\frac{r^{g}}{(r+k)^{g}} \bfv_{g}(r, k).$$ Here, $SU_{X}(r)$ denotes the moduli space of bundles with trivial determinant, while $\mathcal L$ is the unique ample generator of the Picard group. In fact, $\mathcal L$ is the restriction of (an arbitrary) theta bundle $\Theta_{\kappa}$ to $SU_{X}(r)\hookrightarrow U_{X}(r,0)$. 

\subsubsection{} The following genus $g$ symbol is a generalization of Jordan's totient. For any integer $h\geq 2$, we decompose $$h=p_1^{a_1}\ldots p_n^{a_n}$$ into powers of primes. We set \begin{equation}\label{symb}\left\{\frac{\lambda}{h}\right\}_{g}=\begin {cases} 0 & \text {if } p_1^{a_1-1}\ldots p_n^{a_n-1} \text {does not divide } \lambda,\\ {\prod}_{i=1}^{n} \left(\epsilon_{i}-\frac{1}{p_i^{2g}}\right) & \text {otherwise },\end {cases}\end{equation} where $$\epsilon_{i}=\begin {cases} 1 & \text{ if } p_{i}^{a_{i}} |\lambda,\\ 0 & \text {otherwise.} \end {cases}$$ 
If $h=1$, the symbol is always defined to be $1$. A similar symbol appeared in \cite {O} for genus $g=1$.

\subsubsection{} For any degree $0$ line bundle $\xi$ on the Jacobian, we let $\tilde \xi$ be an $r$-th root of $\xi$,  i.e. $\tilde \xi^{\otimes r}=\xi$. The vector bundle $$\bfw_{r, k, \xi}=\bfw_{r,k}\otimes \tilde \xi$$ on the Jacobian only depends on $\xi$. Indeed, the roots $\tilde \xi$ differ by $r$-torsion line bundles $y$; independence of choices follows from the isomorphism $$\bfw_{r,k}\otimes y\cong \bfw_{r,k}$$ explained in Corollary $7.12$ of \cite {mukai} or by equation \eqref{l67} below. 

\subsubsection{} We will prove the following

\begin {theorem}\label{thm2} If $\gcd (r, k)=1$, and $h$ is odd, the Verline bundle splits as \begin{equation}\label{splbfe}\bfe_{hr, hk}\cong \bigoplus_{\xi}\bfw_{r, k, \xi}^{\oplus{\bfm}_{\xi}(r,k)}.\end{equation} For each $h$-torsion line bundle $\xi$ on the Jacobian, having order $\omega$, the multiplicity of the bundle $\bfw_{r, k, \xi}$ in the above decomposition equals \begin{equation}\label{mul1}\bfm_{\xi}(r,k)=\sum_{\delta|h} \frac{1}{(r+k)^{g}\delta^{2g}} \left\{\frac{h/\omega}{h/\delta}\right\}_{g} \bfv_{\frac{h}{\delta}(g-1)+1}({r\delta, k\delta}).\end{equation} \end {theorem}

We note the following features of formula \eqref{mul1}:
\begin {itemize}
\item [(i)] the multiplicities correctly predict the rank of the Verlinde bundles;
\item [(ii)] the pullbacks of both sides by the morphism $hr:A\to A$ are seen to agree \cite {Po};
\item [(iii)] the expression matches the genus $1$ computation of \cite {O};
\item [(iv)] it matches the coprime rank-level splitting \eqref{coprime1}.
\end {itemize}

The case of $h$ even not covered by the theorem involves additional signs which we believe to be related to the Weyl pairing on the Jacobian. 

\subsection{Invariance under Fourier-Mukai} \subsubsection{} The splitting \eqref{splbfe} given by the Theorem is preserved by the Fourier-Mukai transform. To explain this statement, recall first that as a consequence of strange duality, there is a level-rank isomorphism \begin{equation}\label{bfe}\bfe_{r,k}^{\vee}\to \widehat {\bfe}_{k,r}\end{equation} observed by Popa in \cite {Po}. Here, the hat denotes the Fourier-Mukai transform with kernel the normalized Poincar\'e bundle on the Jacobian. 

\subsubsection{} \label{comm} For $k$ and $r$ odd, we will prove in Section \ref{shfm} that the indecomposable factors satisfy the similar duality \begin{equation}\label{bfw}\bfw_{r,k, \xi}^{\vee}\cong \widehat {\bfw_{k, r,\xi}}.\end{equation} Moreover, using \eqref{symmetry1} it is easy to observe the level-rank symmetry of the multiplicities \eqref{mul1}: \begin{equation}\label{bfm}\bfm_{\xi}(r,k)=\bfm_{\xi}(k,r).\end{equation} Therefore, the splitting \eqref{splbfe} is in agreement with Popa's observation \eqref{bfe}. However, these considerations {\it do not} provide an alternative proof of strange duality; the compatibility of the isomorphisms \eqref{splbfe}, \eqref{bfe}, \eqref{bfw} is not entirely obvious.

\subsection {The action of torsion points on the space of generalized theta functions} \subsubsection{} The proof of Theorem \ref{thm2} relies foremost on the computation of the trace of a torsion point in the Jacobian on the space of generalized theta functions. More precisely, for an $h$-torsion line bundle $\alpha$ on the curve $X$, we consider the induced tensoring action $$E\to E\otimes \alpha$$ on the moduli space $SU_X(hr)$. For $h$ odd, the line bundle $\mathcal L^{hk}$ descends to the quotient $SU_X(hr)/\langle \alpha \rangle$ \cite {BLS}, hence it comes equipped with a natural $\alpha$-linearization; see Remark \ref{rem1} in Section \ref{vers} for more details. Therefore, the space of generalized theta functions $$H^0(SU_X(hr), \mathcal L^{hk})$$ decomposes according to the action of $\alpha$. This is the fiberwise splitting \eqref{splbfe} of the Verlinde bundles. 

\subsubsection{} We prove that the trace of $\alpha$ is expressed, somewhat surprisingly, in terms of Verlinde numbers in higher genus and smaller rank:

\begin{theorem}\label{prop1} Let $h$ be odd. For an $h$-torsion line bundle $\alpha$ of order exactly $\delta$ in the Jacobian, we have $$\text {Trace}\left(\alpha, H^{0}(SU_{X}(hr), \mathcal L^{hk})\right)=\frac{r^g}{(k+r)^{g}}{\bfv}_{(g-1){\delta}+1}\left(\frac{hr}{\delta}, \frac{hk}{\delta}\right).$$\end{theorem} 

Note that this result was proved for $g=1$ in \cite {O}. When $\delta=h$, the corresponding formula was obtained by Beauville \cite {B}. Our proof of Theorem \ref{prop1} has the same starting point as Beauville's calculation, but the details are more involved. Related computations for the particular case $h=2$ and $r=1$, not covered by our result, can be found in \cite {AM}.

\subsubsection {}We explain in the last Section that Theorem \ref{prop1} leads to the Verlinde formula for non simply connected groups covered by $SL_{r}$. For the group $PGL_{r}$ with $r$ prime the corresponding formula is due to Beauville \cite {B}. Here, we consider the case of arbitrary $r$. Note that the case of arbitrary non simply connected groups was fully worked out by Alexeev, Meinrenken and Woodward by symplectic methods \cite {AMW}.

\subsection {The semihomogeneous Wirtinger duality} \label{classical} Our study of the semihomogeneous bundles $\bfw_{a,b}$ allows for a generalization of the classical Wirtinger duality. The latter provides an identification of the space of level $2$ theta functions with its dual: \begin{equation}\label{cw}W:H^{0}(A, \Theta^{\otimes 2})^{\vee}\to H^{0}(A, \Theta^{\otimes 2})\end{equation} 

\subsubsection{}\label{classical1} The isomorphism \eqref{cw} is constructed geometrically as follows. Consider the isogeny $$\mu:A\times A\to A\times A, (x,y)\to (x+y, x-y).$$ Using the see-saw theorem, one shows that $$\mu^{\star}(\Theta\boxtimes \Theta)\cong \Theta^{\otimes 2}\boxtimes \Theta^{\otimes 2}.$$ The unique pullback section of $\mu^{\star}(\Theta\boxtimes \Theta)$, namely $$\theta(x+y)\theta(x-y),$$  gives by K\"{u}nneth decomposition an element in $H^{0}(A, \Theta^{\otimes 2})\otimes H^{0}(A, \Theta^{\otimes 2})$ inducing the isomorphism $W$  \cite {mumford}. The explicit K\"{u}nneth decomposition of the above section is tantamount to the addition formulas for theta functions.

 \subsubsection{} \label{wdi}To explain the higher rank duality, we modify the construction of Subsection \ref{classical1} as follows. For odd coprime positive integers $a$ and $b$, consider the morphism $$\mu:A\times A\to A\times A, \,\,\, (x,y)\to (ax+by, x-y).$$ We will prove in Proposition \ref{lemma1} that $$\mu^{\star}(\bfw_{ab,1}\boxtimes \bfw_{1,1})\cong \bfw_{b,a+b}\boxtimes \bfw_{a, a+b}.$$ It is not hard to show that the vector bundle $$\bfw_{ab,1}\boxtimes \bfw_{1,1}$$ has a unique holomorphic section. Consequently, the K\"{u}nneth decomposition of the pullback section induces a morphism \begin{equation}\label{gw}W:H^{0}(A,\bfw_{a,a+b})^{\vee}\to H^{0}(A, \bfw_{b,a+b}).\end{equation} We will prove

\begin {theorem} \label{semih}The morphism $W$ is an isomorphism.
\end {theorem}

Note that the classical Wirtinger duality \eqref{cw} is obtained when $a=b=1$. For this reason, the morphism $W$ of \eqref{gw} will be termed the {\it semihomogeneous Wirtinger duality}. A more general version of Theorem \ref{semih} will be proved in Section \ref{variants}.

\subsubsection{} Theorem \ref{semih} is related to the {\it Generalized Wirtinger Duality} which was the main step of the proof of strange duality given in \cite {MO}. To illustrate this point, assume for simplicity that $\gcd(r, k)=1$. There is an isomorphism \begin{equation}\label{gwd}{\mathsf D}: H^{0}(U_{X}(r,0), \Theta^{k}\otimes {\det}^{\star}\Theta)^{\vee}\to H^{0}(U_{X}(k,0), \Theta^{r}\otimes {\det}^{\star}\Theta).\end{equation} Pushing forward by the determinant, we thus obtain an isomorphism $$H^{0}(A, \bfe_{r,k}\otimes \Theta)^{\vee}\to H^{0}(A, \bfe_{k,r}\otimes \Theta).$$ This is in agreement with the splitting \eqref{coprime1} and the isomorphism $$H^{0}(A, \bfw_{r,k}\otimes \Theta)^{\vee}=H^{0}(A, \bfw_{r, k+r})^{\vee}\to H^{0}(A, \bfw_{k, k+r})=H^{0}(A, \bfw_{k,r}\otimes \Theta)$$ explained by Theorem \ref{semih}. However, this argument does not immediately provide an alternative proof of \eqref{gwd}, for the reasons mentioned in Subsection \ref{comm}.

\subsection{Outline} We begin by proving a number of results about semihomogeneous vector bundles in the section immediately following the Introduction. Theorem \ref{semih} will be derived in Section \ref{shwd}. The central part of the paper is the calculation leading to Theorem \ref{prop1}; this will be presented in Section \ref{tracecalc}. We will show in Section \ref{vers} how Theorem \ref{thm2} follows from this calculation. The last section explains the connection between Theorem \ref{prop1} and the Verlinde formula for non simply connected groups. 

\subsection{Acknowledgments} The author gratefully acknowledges conversations with Alina Marian. This work was partially supported by the NSF grant DMS-0852468. 
 
\section {Generalities on semihomogeneous bundles} 
In this section, we collect several general facts about semihomogenous bundles, following Mukai's foundational work \cite{mukai}. 

\subsection {The bundles $\bfw_{a,b}$.} To begin, let $(A, \Theta)$ be a principally polarized complex abelian variety with $\Theta$ symmetric. Let $a, b$ be odd coprime positive integers. By Theorem $7.11$ and Remark $7.13$ of \cite {mukai}, there exist simple semihomogeneous vector bundles $\bfw_{a,b}$ such that  \begin{equation}\label{numinv} \text {rank }\bfw_{a,b}= a^{g}, \,\, \det \bfw_{a,b}=\Theta^{a^{g-1}b}.\end{equation} Such bundles are not unique; they differ by twists by suitable line bundles on $A$. We will show existence and uniqueness under the symmetry assumption $$(-1)^{\star}\bfw_{a,b}\cong \bfw_{a,b}.$$ 

To this end, let $\bfw$ be any simple, possibly non-symmetric, semihomogeneous bundle as in \eqref{numinv}. Since $(-1)^{\star}\bfw$ is also simple, semihomogeneous and satisfies \eqref{numinv}, Theorem $7.11$ in \cite {mukai} implies that $$(-1)^{\star}\bfw\cong \bfw\otimes M,$$ for some degree $0$ line bundle $M$ on $A$. By comparing determinants, we obtain $$M^{\otimes a^{g}}=0.$$ Since $a$ is odd, we can pick an $a^{g}$-torsion line bundle $L$ on $A$ such that $L^{\otimes 2}=M$. Then, $$\bfw'=\bfw\otimes L$$ clearly satisfies \eqref{numinv}. Moreover, $\bfw'$ is symmetric: $$(-1)^{\star}\bfw'\cong(-1)^{\star}\bfw\otimes L^{-1}\cong\bfw\otimes M\otimes L^{-1}\cong\bfw\otimes L\cong\bfw',$$ thus proving existence. 

To show uniqueness, observe that any two vector bundles $\bfw_{1}$, $\bfw_{2}$ with numerical invariants \eqref{numinv} must satisfy $$\bfw_{1}\cong\bfw_{2}\otimes M$$ for an $a^{g}$-torsion line bundle $M$, again by Theorem $7.11$ \cite {mukai}. Symmetry implies that $$\bfw_{1}\cong(-1)^{\star}\bfw_{1}\cong(-1)^{\star}\bfw_{2}\otimes M^{-1}\cong\bfw_{2}\otimes M^{-1}\cong\bfw_{1}\otimes M^{-2}.$$ By \cite {mukai}, Corollary $7.12$, or by equation \eqref{l67} below, the group $$\Sigma(\bfw_1)=\{\text{line bundle } P \text { such that } \bfw_1\otimes P\cong \bfw_1\}$$ is isomorphic to the group of $a$-torsion line bundles on $A$. We conclude that $$M^{\otimes 2a}=M^{\otimes a^{g}}=0.$$ Since $a$ is odd, it follows that $M^{\otimes a}=0$. Hence $M\in \Sigma(\bfw_1)$, and therefore $$\bfw_{1}\cong\bfw_{2}\otimes M\cong\bfw_{2},$$ thus proving uniqueness.

\subsection {Translation invariance and Heisenberg actions}\label{complred} Much like line bundles, semihomogeneous vector bundles have associated Theta groups. We will make this statement precise now.

\subsubsection {} \label{hact}First, recall that for any $x\in A$, there is a line bundle $y$ on $A$ such that $$t_{x}^{\star}\bfw_{a,b}=\bfw_{a,b}\otimes y.$$ In fact, by Lemma $6.7$ of \cite {mukai}, we may assume that \begin{equation}\label{l67}x=a\alpha \,\text { and }y=b\alpha\text { for } \alpha \in A.\end{equation} Here, $y$ is viewed as a point of $A$ by means of the principal polarization $\Theta$ which identifies the abelian variety with its dual. It follows that if $x$ is a $b$-torsion point, then we may pick $\alpha$ to be a $b$-torsion point as well. Thus \begin{equation}\label{translation}t_{x}^{\star}\bfw_{a,b}\cong \bfw_{a,b}.\end{equation} Conversely, if \eqref{translation} is satisfied, then $x$ is a $b$-torsion point, as explained in Corollary $7.12$ of \cite {mukai}.

We let $A[b]$ denote the group of $b$-torsion points on $A$. Let $\bfh$ be the group of pairs $(x,\Phi)$, where $x\in A[b]$ and $$\Phi:\bfw_{a,b}\to t_x^{\star}\bfw_{a,b}$$ is an isomorphism. The group law on $\bfh$ is given by $$(x, \Phi)\cdot (y, \Psi) =(x+y, t_y^{\star}\Phi\circ \Psi).$$ Note that there is an extension $$0\to \mathbb C^{\star}\to \bfh\to A[b]\to 0.$$ The group $\bfh$ is a nondegenerate Heisenberg group, {\it i.e.} $\mathbb C^{\star}$ is the center of $\bfh$. Nondegeneracy follows by the usual arguments, see for instance Theorem $1$ in \cite{Mu} or the discussion in Section $2$ of \cite {Um}.

\subsubsection{} Crucially for our arguments, any nondegenerate Heisenberg group has a unique irreducible Schr\"odinger representation $\bfs$ for which the center $\mathbb C^{\star}$ acts with the natural character. Moreover, we have $$\dim \bfs = b^{g}.$$ Any other representation for which the center acts by the natural character is direct sum of copies of $\bfs$, \cite {Mu}.

The group $\bfh$ introduced in Subsection \ref{hact} acts on the vector bundle $\bfw_{a,b}$, as well as on its holomorphic sections: $$H^0(A, \bfw_{a,b})\times \bfh \ni (s, (x, \Phi))\mapsto t_{-x}^{\star} \Phi(s)\in H^0(A, \bfw_{a,b}).$$ Furthermore, the center of $\bfh$ acts by homotheties. Comparing dimensions via \eqref{higher1} below, it follows that $H^0(A, \bfw_{a,b})$ is a geometric realization of the Schr\"{o}dinger representation $\bfs$. 

\subsubsection{} When $a=1$ the above remarks specialize as follows. For any integer $m$, let $\bfh[m]$ be the Heisenberg group of the line bundle $\Theta^{m}$ obtained as an extension $$0\to \mathbb C^{\star}\to \bfh[m]\to A[m]\to 0.$$ A variation of this exact sequence will be useful throughout the paper. When $m$ is odd, the assignment $$\eta\to \eta^m$$ is an endomorphism of $\bfh[m]$ whose image lies in the center $\mathbb C^{\star}$. We define $\widetilde \bfh[m]$ to be the kernel of this endomorphism sitting in the exact sequence $$0\to \mu_m\to \widetilde \bfh[m]\to A[m]\to 0.$$ When $m$ is even, one may repeat this construction starting with the homomorphism $\eta\to \eta^{2m}.$

The irreducible representations of the finite group $\widetilde \bfh[m]$ have been worked out in \cite {S}. In particular, if $\gcd(m,n)=1$, there is a unique irreducible Schr\"{o}dinger representation $\bfs_{m,n}$ of dimension $$\dim \bfs_{m,n}=m^{g}$$ such that the center acts by the character $\alpha\to \alpha^{n}.$  Since the center $\mathbb C^{\star}$ and the finite subgroup $\widetilde \bfh[m]$ generate the group $\bfh[m],$ we immediately see that $\bfs_{m,n}$ becomes the unique irreducible representation of $\bfh[m]$ such that the center acts with weight $n$. 

\subsection {Splittings under isogenies} \label{splittingsisog} It is a general fact established in \cite {mukai} that pullbacks of semihomogeneous bundles under suitable isogenies split. This can be checked directly for the bundles $\bfw_{a,b}$ studied here.

\subsubsection{} Consider the isogeny $a:A\to A.$ We claim that \begin{equation}\label{splittings}a^{\star}\bfw_{a,b}=\oplus_{i=1}^{a^{g}} \Theta^{ab}.\end{equation} This is a consequence of the construction of the bundles $\bfw_{a,b}$. Indeed, for any semihomogeneous bundle $\bfw$, set $$\Phi(\bfw)=\{(x, y): t_{x}^{\star}\bfw=\bfw \otimes y\}\subset A\times \widehat A.$$ By Lemma $3.6$ of \cite {mukai}, the pullback under the projection $$p:\Phi (\bfw)\to A$$ splits as $$p^{\star}\bfw\cong \bigoplus L,$$ for some line bundle $L$. In our case, the assignment $$A\to \Phi (\bfw_{a,b}), \alpha \mapsto (a\alpha, b\alpha)$$ is an isomorphism cf. Proposition $7.7$ of \cite {mukai}; see also \eqref{l67}. The projection $p:A\cong\Phi(\bfw_{a,b})\to A$ becomes multiplication by $a$. Thus we obtain $$a^{\star}\bfw_{a,b}\cong\oplus_{1}^{a^{g}}L.$$ Taking determinants, we see that $$L\cong\Theta^{ab}\otimes M$$ where $M$ is an $a^{g}$-torsion point. However, observing that both $\bfw_{a,b}$ and $\Theta$ are symmetric, we derive that $L$ and $M$ are symmetric. Since $a$ is odd, this implies that $M$ is trivial, proving the claim.

\subsubsection{} We sharpen equation \eqref{splittings} by including Heisenberg group actions. The isogeny $$a:A\to A$$ has as kernel the group $A[a].$ Therefore the bundle $a^{\star}\bfw_{a,b}$ carries a natural action of the Heisenberg group $\bfh[a]$, such that the center acts trivially in the fibers over the identity. On the other hand, the bundle $\Theta^{ab}=(\Theta^{\otimes a})^{\otimes b}$ carries a linearization of $\bfh[a]$ such that the center acts as $\alpha\to \alpha^{b}$ in the fibers over the identity. We claim that \begin{equation}\label{heisenberg}a^{\star}\bfw_{a,b}\cong \Theta^{ab}\otimes \bfs_{a,-b}\end{equation} as representations of $\bfh[a]$. Indeed, $$H^{0}(A, a^{\star}\bfw_{a,b}\otimes \Theta^{-ab})$$ becomes a representation of $\bfh[a]$ of dimension $a^{g}$ such that the center acts as $\alpha\to \alpha^{-b}.$ The claim follows by the uniqueness of the Schr\"{o}dinger representation.

Finally, note that equation \eqref{heisenberg} may be taken as the definition of $\bfw_{a,b}$ when either $a$ or $b$ is even. 

\subsection {Higher cohomology} In general, the cohomology of semihomogeneous bundles with non-zero Euler characteristic vanishes in all but one degree \cite {Um}. This is easily seen for the bundles $\bfw_{a,b},$ as a consequence of the discussion in Subsection \ref{splittingsisog}. We have $$h^{i}(A, \bfw_{a,b})=0$$ if $b>0$, $i>0$; this statement is obvious after pullback by the \'etale cover $a$. By Remark $7.13$ in \cite{mukai}, we have $$\chi(A, {\bfw}_{a,b})=b^{g}.$$ Thus \begin{equation}\label{higher1}h^{0}(A, \bfw_{a,b})=b^{g}.\end{equation} In particular, the vector bundle $\bfw_{a,1}$ has a unique holomorphic section.

\section{The semihomogeneous Wirtinger Duality}
\label{shwd}

In this section we will prove Theorem \ref{semih}. As a byproduct of the proof, we show the invariance of the semihomogeneous bundles under the Fourier-Mukai transform.

\subsection {The duality and its proof} \label{dpfa} The semihomogeneous Wirtinger map was set up in Subsection \ref{wdi}. For odd coprime integers $a,b$, recall the isogeny $$\mu:A\times A\to A\times A, \,\,\, (x,y)\to (ax+by, x-y).$$
\subsubsection{} \label{dpf} We claim that \begin {proposition}\label{lemma1} We have 
$$\mu^{\star}(\bfw_{ab,1}\boxtimes \bfw_{1,1})\cong\bfw_{b,a+b}\boxtimes \bfw_{a, a+b}.$$
\end {proposition}
{\it Proof:} The vector bundle $$\bfw_{b,a+b}\boxtimes \bfw_{a,a+b}$$ on the abelian variety $A\times A$ is semihomogeneous. Its rank equals $(ab)^{g}$, while the determinant is $$\Theta^{a^{g}b^{g-1}(a+b)}\boxtimes \Theta^{b^{g}a^{g-1}(a+b)}=\left(\Theta^{a(a+b)}\boxtimes\Theta^{b(a+b)}\right)^{(ab)^{g-1}}.$$ Similarly, $$\mu^{\star}(\bfw_{ab,1}\boxtimes \bfw_{1,1})$$ is semihomogeneous of rank $(ab)^{g}$ and determinant $$\mu^{\star}\left(\Theta^{(ab)^{g-1}}\boxtimes \Theta^{(ab)^{g}}\right)=\mu^{\star}(\Theta\boxtimes\Theta^{ab})^{(ab)^{g-1}}.$$ We claim that $$\mu^{\star}(\Theta\boxtimes\Theta^{ab})\cong \Theta^{b(a+b)}\boxtimes\Theta^{a(a+b)}.$$ This isomorphism is a consequence of the see-saw theorem. Indeed, the restriction of the pullback $$\mu^{\star}(\Theta\boxtimes\Theta^{ab})$$ to $A\times \{y\}$, for $y\in A$, equals $$a^{\star}t_{by}^{\star}\Theta\otimes t_{-y}^{\star}\Theta^{ab}=t_{by/a}^{\star}a^{\star}\Theta\otimes t_{-y}^{\star}\Theta^{ab}=\left(t_{by/a}^{\star}\Theta^{a}\otimes t_{-y}^{\star}\Theta^{b}\right)^{a}=\left(\Theta^{a+b}\right)^{a}=\Theta^{a(a+b)}.$$ In the third equation, we made use of the fact that $\Theta$ is symmetric to conclude that $a^{\star}\Theta=\Theta^{a^{2}}.$ For the last equality, we used the theorem of the square. The computation of the restriction to $\{x\}\times A$ is similar. Therefore, the ranks and determinants of the two semihomogeneous bundles $$\mu^{\star}(\bfw_{ab,1}\boxtimes \bfw_{1,1}) \text{ and } \bfw_{b,a+b}\boxtimes \bfw_{a,a+b}$$ agree. 

We will prove that the two bundles are isomorphic by using the structure theory of semihomogeneous bundles developed in Theorem $7.11$ and Proposittion $6.18$ of \cite {mukai}. Let $\delta={\det}/{\text{rank}}$ denote the common slope of the two vector bundles above. The category of semihomogeneous vector bundles of slope $\delta$ contains a (nonunique) simple bundle $E$. Any semihomogeneous bundle $F$ of slope $\delta$, not necessarily simple, can be expressed as direct sum $$F=\bigoplus F_i,$$ where $F_i$ are semihomogeneous. Moreover, the $F_i$'s admit filtrations $$G_1^i\subset G_2^i\subset \ldots \subset G_n^i$$ whose successive quotients are of the form $$G_j^i/G^i_{j-1}\cong E\otimes M_i,$$ for suitable degree $0$ line bundles $M_i$. 

In our case, we may take $$E=\bfw_{b,a+b}\boxtimes \bfw_{a,a+b},$$ as this bundle is clearly simple. Applying the above remarks to $$F=\mu^{\star}(\bfw_{ab,1}\boxtimes \bfw_{1,1})$$ and comparing the ranks of $E$ and $F$, we conclude that $$F=E\otimes M,$$ for some degree $0$ line bundle $M$ on $A\times A$. Taking determinants we obtain \begin {equation}\label{eqnm}M^{\otimes (ab)^g}=0.\end{equation} Furthermore, $E$ and $F$ are clearly symmetric, hence $$E\cong (-1)^\star E\cong (-1)^{\star}F\otimes M\cong F\otimes M\cong E\otimes M^{\otimes 2}.$$ Therefore, $M^{\otimes 2}$ belongs to the group $$\Sigma(E)=\{P \text { line bundle on } A\times A\text { such that }E\otimes P\cong E\}.$$ Proposition $7.1$ of \cite {mukai} guarantees that $\Sigma(E)$ has $(\text {rank } E)^2=(ab)^{2g}$ elements. In fact, these elements are seen to correspond to the line bundles of order $ab$ on $A$, via the assignment $$\widehat A[ab]\to \Sigma(E), \,\,L\mapsto P= L^{\otimes a}\boxtimes L^{\otimes b}.$$ Indeed, it suffices to observe that, as $(a,b)=1$, the above assignment is injective, if well defined. To check this last point, we first note that $L^{\otimes a}$ has order $b$ on $A$. Applying Corollary $7.12$ of \cite {mukai} once again, we derive $$\bfw_{b,a+b}\otimes L^{\otimes a}=\bfw_{b,a+b}.$$ Similarly, $$\bfw_{a,a+b}\otimes L^{\otimes b}=\bfw_{a,a+b}.$$ Hence, tensoring by $P=L^{\otimes b}\boxtimes L^{\otimes a}$ preserves the bundle $E$, as claimed. We conclude that $$M^{\otimes 2}=L^{\otimes a}\boxtimes L^{\otimes b},$$ for some $ab$-torsion line bundle $L$ on $A$. We claim that $M$ has the same form. Indeed, set $$\tilde L=L^{\otimes \frac{(ab)^g+1}{2}}.$$ By invoking \eqref{eqnm} we obtain $$M=M^{(ab)^g+1}=\tilde L^{\otimes a}\boxtimes \tilde L^{\otimes b}.$$ Therefore $M\in \Sigma(E),$ hence $$E\cong E\otimes M\cong F.$$ The Proposition is proved.

\qed

\vskip.1in
\subsubsection{} As explained in Subsection \ref{wdi}, the unique tensor product section of $\bfw_{ab,1}\boxtimes \bfw_{1,1}$ induces by pullback and K\"unneth decomposition the semihomogeneous Wirtinger duality map \begin{equation}\label{wdm}W:H^{0}(A, \bfw_{a, a+b})^{\vee}\to H^{0}(A, \bfw_{b, a+b}).\end{equation} We show that $W$ is an isomorphism. 
\vskip.1in

{\it Proof of Theorem \ref{semih}:} First, note that by \eqref{higher1} the two sides of \eqref{wdm} have the same dimension, namely $(a+b)^{g}$. 

The isogeny $\mu$ is invariant by the translation action of $A[a+b]$ on both factors. Therefore, the bundle $$\mu^{\star}(\bfw_{ab,1}\boxtimes \bfw_{1,1})\cong \bfw_{b,a+b}\boxtimes \bfw_{a,a+b}$$ carries a natural linearization of $A[a+b]$. In particular, for any $x\in A[a+b]$, we have natural identifications $$\chi_x:\bfw_{b,a+b}\boxtimes \bfw_{a,a+b}\to t_{x}^{\star}(\bfw_{b,a+b}\boxtimes \bfw_{a,a+b})$$ which are compatible with the group structure of $A[a+b]$.

By \eqref{translation}, the bundles $\bfw_{a,a+b}$ and $\bfw_{b,a+b}$ are invariant under the translation action of $A[a+b]$. Let $\bfh$ denote the Theta group of $\bfw_{a,a+b}$ sitting in an exact sequence $$0\to \mathbb C^{\star}\to \bfh\to A[a+b]\to 0.$$ We explained in Subsection \ref{hact} that $\bfw_{a,a+b}$ admits an $\bfh$-linearization such that the center $\mathbb C^{\star}$ acts with weight $1$ in the fiber over the identity. We claim that we can pick an $\bfh$-linearization of $\bfw_{b,a+b}$ such that center acts with opposite weight $-1$ in the fiber over the identity. Indeed, any $h\in \bfh$ gives a torsion point $x\in A[a+b]$ and an isomorphism $$\Phi:\bfw_{a,a+b}\to t_{x}^{\star}\bfw_{a,a+b}.$$ In order to define the action of $h$ on $\bfw_{b,a+b}$ we specify an isomorphism $$\Psi:\bfw_{b, a+b}\to t_{x}^{\star}\bfw_{b,a+b}$$ such that the product $$\Psi\boxtimes \Phi:\bfw_{b,a+b}\boxtimes \bfw_{a,a+b}\to t_{x}^{\star}\bfw_{b,a+b}\boxtimes t_{x}^{\star}\bfw_{a,a+b}$$ equals $\chi_x.$ To construct $\Psi$, first pick an arbitrary isomorphism $\Psi_{0}$, and observe that $$\chi_x^{-1}\circ \Psi_{0}\boxtimes \Phi:\bfw_{b,a+b}\boxtimes \bfw_{a,a+b}\to \bfw_{b,a+b}\boxtimes \bfw_{a,a+b}$$ is an automorphism of a simple bundle, hence it is a multiple of the identity. It suffices to define $\Psi$ as the appropriate rescaling of $\Psi_{0}$. The construction is independent of choices.

A dimension count shows that the two vector spaces $$H^{0}(A, \bfw_{a,a+b})^{\vee},\,\, H^{0}(A, \bfw_{b,a+b})$$ are the irreducible Schr\"{o}dinger representation $\bfs_{a+b, -1}$ of $\bfh$. The Heisenberg invariance of the morphism $W$ gives the statement. \qed 

\subsection {Variants} \label{variants}A slightly more general variant of Theorem \ref{semih} will be given here. Consider the matrix \[ {\bf M}=\left( \begin{array}{cc}
a & b \\
c & -d \end{array} \right)\] where $a,b,c,d$ are odd positive integers such that $(ac,bd)=1$. Let $$\delta=-\det {\bf M}=ad+bc.$$  The case we considered in Section \ref{dpfa} corresponds to $c=d=1$. 

Consider the morphism $$\mu: A \times A \to A\times A,$$ \[\left[\begin{array}{c} x \\y\end{array}\right]\to {\bf M}\left[\begin{array}{c} x\\ y\end {array}\right]=\left[\begin{array}{c}ax+by\\ cx-dy\end{array}\right].\]Then 
\begin{lemma} There is an isomorphism $$\mu^{\star}(\bfw_{ab,1}\boxtimes \bfw_{cd,1})\cong \bfw_{bd, \delta} \boxtimes \bfw_{ac, \delta}.$$\end {lemma} 

The following can be proved in a completely similar fashion as Theorem \ref{semih}:
\begin {theorem}The unique section of $$\bfw_{ab,1}\boxtimes \bfw_{cd,1}$$ induces by pullback and K\"unneth decomposition an isomorphism $$W: H^{0}(A, \bfw_{bd, \delta})^{\vee}\to H^{0}(A, \bfw_{ac, \delta}).$$\end {theorem}

\subsection {The semihomogeneous bundles and Fourier-Mukai}\label{shfm}

The discussion in Subsection \ref{dpf} also gives the isomorphism \eqref{bfw} of the Introduction: $$\widehat {\bfw_{a,b, \xi}}\cong \bfw_{b,a,\xi}^{\vee},$$ for any degree $0$ line bundle $\xi$. In fact, it suffices to prove the case of trivial $\xi$: \begin{equation}\label{l12}\widehat{\bfw}_{a,b}\cong \bfw_{b,a}^{\vee}.\end{equation} The general statement follows immediately from here. Indeed, pick a root $\eta^{ab}=\xi$. Using the propreties of Fourier-Mukai and equation \eqref{l67}, we compute 
$$\widehat {\bfw_{a,b, \xi}}\cong \widehat {\bfw_{a,b}\otimes \eta^{b}}\cong t_{b\eta}^{\star} \widehat{\bfw}_{a,b}\cong t_{b\eta}^{\star}\bfw_{b,a}^{\vee}\cong \bfw_{b,a}^{\vee}\otimes \eta^{-a}\cong\bfw_{b,a, \xi}^{\vee}.$$ 

To prove \eqref{l12}, consider the morphism $$f:A\times A\to A, (x,y)\mapsto ax+by.$$ Let $\mathcal P$ be the Poincare bundle on the product $A\times A$, normalized in the usual way.
We have the following

\begin {lemma} \begin{itemize}\item[(i)]There is an isomorphism $$f^{\star}\bfw_{ab,1}=\left(\bfw_{b,a}\boxtimes \bfw_{a,b}\right)\otimes \mathcal P.$$
 \item [(ii)] There is an isomorphism $$\bfw_{b,a}^{\vee}\cong \widehat \bfw_{a,b},$$ induced by the unique pullback section of $f^{\star}\bfw_{ab,1}$ via the identification $$H^0(A\times A, f^{\star}\bfw_{ab,1})\cong H^0(A\times A, \bfw_{b,a}\boxtimes \bfw_{a,b}\otimes \mathcal P)\cong H^0(A, \bfw_{b,a}\boxtimes \widehat \bfw_{a,b}).$$ 
 \end {itemize}
 \end {lemma} 

{\it Proof.} Let $$g:A\times A\to A, (x,y)\to x-y.$$ From the see-saw principle we conclude that $$g^{\star}\Theta=\left(\Theta\boxtimes \Theta\right)\otimes \mathcal P^{-1}.$$ Item (i) follows from Proposition \ref{lemma1}: $$f^{\star}\bfw_{ab,1}=\mu^{\star}(\bfw_{ab,1}\otimes \Theta)\otimes g^{\star}\Theta^{-1}=(\bfw_{b,a+b}\boxtimes \bfw_{a,a+b})\otimes (\Theta^{-1}\boxtimes \Theta^{-1})\otimes \mathcal P$$ $$=(\bfw_{b,a}\boxtimes \bfw_{a,b})\otimes \mathcal P.$$

As for (ii), note that the pullback section of $f^{\star}\bfw_{ab,1}$ gives a nonzero map \begin{equation}\label{bfwba} \bfw_{b,a}^{\vee}\to \widehat \bfw_{a,b}.\end{equation} We check that both sides are stable and have the same numerical invariants. This implies that \eqref{bfwba} is an isomorphism.

To this end, we first note that both bundles in \eqref{bfwba} are simple and semihomogeneous, hence stable with respect to any polarization \cite {mukai}. To prove the semihomogeneity of $\widehat {\bfw}_{a,b}$, let $x\in A$ and write $x=b\alpha$ for some $\alpha\in A$. Set $y=-a\alpha$. Using the propreties of Fourier-Mukai and \eqref{l67}, we compute $$t_{x}^{\star}\widehat {\bfw}_{a,b}\cong \widehat {\bfw_{a,b}\otimes x}\cong \widehat{\bfw_{a,b}\otimes b\alpha}\cong \widehat {t_{a\alpha}^{\star}\bfw_{a,b}}\cong \widehat {t_{-y}^{\star}\bfw_{a,b}}\cong \widehat{\bfw}_{a,b}\otimes y.$$ To match the numerical invariants of $\bfw_{b,a}^{\vee}$ and $\widehat {\bfw}_{a,b}$, note first that \eqref{splittings} gives $$\text {ch} (\bfw_{a,b})=a^g \exp \left(\frac{b}{a}\Theta \right).$$ By Grothendieck-Riemann-Roch, we compute $$\text{ch} (\widehat{\bfw}_{a,b})=\text{ch} p_{!} (\bfw_{a,b}\otimes \mathcal P)=p_{!}(\text {ch }\bfw_{a,b} \exp (c_1(\mathcal P)))=a^g p_{!}\left( \exp \left(\frac{a}{b}\Theta+c_1(\mathcal P)\right)\right)$$ $$=b^{g} \exp \left(-\frac{b}{a}\Theta\right)=\text{ch} (\bfw_{b,a}^{\vee}).$$ This completes the proof. \qed

\section{The splitting of Verlinde bundles}
\label{vers}

In this section we determine the splitting of the Verlinde bundles, thus proving Theorem \ref{thm2}. The argument relies on the trace calculation of Theorem \ref{prop1} to be explained in Section \ref{tracecalc}. 

\subsubsection{}\label{401} Keeping the same notations as in the Introduction, we begin by considering the following warm-up case:

\begin {lemma}\label{coprime} Let $(r,k)=1$ be coprime positive integers. We have $$\bfe_{r,k}=\bigoplus \bfw_{r,k}.$$
\end {lemma}

{\it Proof:} It is enough to check this equality $A[r]$-equivariantly, after pullback by the isogeny $r:A\to A$. Equivalently, invoking \eqref{heisenberg}, we need to prove that $$r^{\star}\bfe_{r,k}=\bigoplus \bfs_{r, -k}\otimes \Theta^{rk}.$$ We use the fiber diagram 
\begin{center}$\xymatrix{SU_{X}(r)\times A \ar[r]^{t} \ar[d]^{p} & U_{X}(r, 0) \ar[d]^{\det} \\ A\ar[r]^{r} &  A}$\end{center} where $t$ denotes tensor product $(E, L)\mapsto E\otimes L.$ The splitting \begin{equation}\label{pull3}t^{\star}\Theta=\mathcal L\boxtimes \Theta^{r}\end{equation} is well-known. Thus $$r^{\star}\bfe_{r,k}=p_{\star}(\mathcal L^{k}\boxtimes \Theta^{rk})=H^{0}(SU(r), \mathcal L^{k})\otimes \Theta^{rk}.$$ 

The group $\bfh[r]$ acts naturally on $\Theta^{r}$ and on $t^{\star}\Theta$. Much like in the proof of Theorem \ref{semih}, the isomorphism \eqref{pull3} gives an action of $\bfh[r]$ on $\mathcal L$, such that the center acts as $\alpha\to \alpha^{-1}$ in the fiber over the identity. Thus, the Heisenberg acts on $H^{0}(SU(r), \mathcal L^{k}),$ such that the center acts with weight $-k$. Since $(r,k)=1$, the representation $H^0(SU(r), \mathcal L^k)$ splits as direct sum of $\bfs_{r,-k}$, by uniqueness of the Schr\"{o}dinger representation and complete reducibility (see Section \ref{complred}). This concludes the proof. \qed

\vskip.1in

\subsubsection{} A variation of the above argument gives the general case.  \vskip.1in

{\it Proof of Theorem \ref{thm2}:} Let us write $hr$ for the rank and $hk$ for the level. It suffices to check that $A[hr]$-equivariantly we have an isomorphism of vector bundles \begin{equation}\label{pullhr}(hr)^{\star}\bfe_{hr, hk}\cong \bigoplus_{\xi} (hr)^\star\bfw_{r,k, \xi}^{\oplus \bfm_{\xi}},\end{equation} where $\xi$ ranges over the characters of $A[h]$, or equivalently over the $h$-torsion line bundles on the Jacobian. We will confirm the multiplicities $$\bfm_{\xi}=\sum_{\delta|h} \frac{1}{(r+k)^{g}\delta^{2g}} \left\{\frac{h/\omega}{h/\delta}\right\}_{g} \bfv_{\frac{h}{\delta}(g-1)+1}(r{\delta},k{\delta}),$$ for characters $\xi$ of order exactly $\omega$.

Recall that we denoted by $\widetilde \bfh[{hr}]$ the finite Heisenberg group of the line bundle $\Theta^{hr}$ on the Jacobian $$0\to \mu_{hr}\to 
\widetilde \bfh[{hr}]\to A[hr]\to 0.$$ (For $r$ even, we need to replace the leftmost group by $\mu_{2hr}.$) As in Lemma \ref{coprime}, we will consider the action of $\widetilde \bfh[hr]$ on the objects in sight, in particular on the pair $(SU_{X}(hr), \mathcal L)$. We have seen in the proof of Lemma $3$ that $$(hr)^{\star}\bfe_{hr,hk}\cong H^{0}(SU_{X}(hr), \mathcal L^{hk})\otimes (\Theta^{hr})^{hk}.$$ Furthermore, by \eqref{heisenberg}, the pullback $(hr)^{\star}\bfw_{r,k}$ splits as sums of pluri-theta bundles. In fact, $\widetilde \bfh[hr]$-equivariantly, we have $$(hr)^{\star} \bfw_{r,k}\cong \left(\Theta^{hr}\right)^{hk}\otimes R,$$ where $R$ is a representation of $\widetilde \bfh[{hr}]$ of dimension $r^g$ such that the center acts with weight $-h k.$ Therefore, to prove \eqref{pullhr} it suffices to check that $\widetilde \bfh[hr]$-equivariantly we have \begin{equation}\label{decomp} H^{0}(SU_{X}(hr), \mathcal L^{hk})=R\otimes \bigoplus_{\xi} \xi^{\bfm_{\xi}}.\end{equation} 

To establish \eqref{decomp}, we will make use of the following crucial observation. Consider the natural morphism of Theta groups $$\iota:\widetilde \bfh[h]\to \widetilde \bfh[hr].$$ On the centers, $\iota$ restricts to $$\alpha\to \alpha^r.$$ Consider two irreducible representations $R_1, R_2$ of $\widetilde \bfh[hr]$ with the same central weight $-hk$. Via $\iota$, $R_1$ and $R_2$ become representations of $\widetilde \bfh[h]$ such that the center acts trivially. Hence $R_1$ and $R_2$ are representations of the group $A[h]$. It follows from equation ($50$) in \cite {O} that
$$R_1\cong R_2 \text { as representations of } \widetilde \bfh[hr] \text { iff }  R_1\cong R_2 \text { as representations of } A[h].$$ In the light of this observation, it suffices to show that \eqref{decomp} holds as an equality of $A[h]$-modules. 

We describe the action of $A[h]$ on the left hand side of \eqref{decomp}. First, recall from \eqref{pull3} the $A[h]$-equivariant identification \begin{equation}\label{pull4}t^{\star}\Theta^h\cong \mathcal L^h\boxtimes \left(\Theta^{hr}\right)^h\cong \mathcal L^h\boxtimes h^{\star}\Theta^r.\end{equation} The last isomorphism uses the fact that $h$ is odd, to conclude that $A[h]$-equivariantly we have $(\Theta^{hr})^h\cong h^{\star}\Theta^r;$ see the the last section of \cite {O}.  As a consequence of \eqref{pull4}, $\mathcal L^h$ acquires an $A[h]$-linearization. 

Regarding the right hand side of \eqref{decomp}, it was explained in Section $4$ of \cite {O} that for $h$ odd, $R$ is the trivial $A[h]$-representation of dimension $r^g$. Therefore, we compute $$\bfm_{\xi}=\frac{1}{r^g h^{2g}}\sum_{\alpha\in A[h]}\xi(\alpha^{-1})\text {Trace} \left(\alpha, H^{0}(SU_{X}(hr), \mathcal L^{hk})\right).$$ Theorem \ref{prop1} shows that $$\text {Trace}\left(\alpha, H^{0}(SU_{X}(hr), \mathcal L^{hk})\right)=\frac{r^g}{(r+k)^{g}}\bfv_{(g-1)\frac{h}{\delta}+1}\left(r{\delta}, k{\delta}\right),$$ for an element $\alpha$ of order $\frac{h}{\delta}$. Then, $$\bfm_{\xi}=\frac{1}{h^{2g}}\sum_{\delta|h} \frac{1}{(r+k)^{g}}\bfv_{(g-1)\frac{h}{\delta}+1}\left(r{\delta}, k{\delta}\right)\left(\sum_{\text{ord}\alpha={h}/{\delta}} \xi(\alpha^{-1})\right)$$ $$=\sum_{\delta|h}\frac{1}{(r+k)^{g}\delta^{2g}}\left\{\frac{h/\omega}{h/\delta}\right\}_{g}\bfv_{(g-1)\frac{h}{\delta}+1}(r\delta, k\delta).$$ In the last line we used (the appropriate extension of)  Lemma $4$ of \cite {O} which gives the value of the sum \begin{equation}\label{sum}\sum_{\text{ord}\alpha={h}/{\delta}}\xi(\alpha^{-1})=\frac{h^{2g}}{\delta^{2g}}\left\{\frac{h/\omega}{{h}/{\delta}}\right\}_{g}.\end{equation} The proof of the Theorem is completed.\qed

\begin{remark} \label{rem1} For $h$ odd, the line bundle $\mathcal L^h$ descends to the quotient $SU_X(hr)/A[h]$, cf. Section $6$ of \cite {BLS}. Therefore, $\mathcal L^h$ admits a natural $A[h]$-linearization such that any element $\alpha\in A[h]$ acts trivially in the fibers of $\mathcal L^h$ over the $\alpha$-fixed points. We claim that this linearization coincides with the one constructed in the proof of Theorem \ref{thm2} via \eqref{pull4}. 

To check this fact, note that since two linearizations differ by a character of $A[h]$,  it suffices to prove that under the isomorphism \eqref{pull4} the action of $\alpha\in A[h]$ in the fibers of $\mathcal L^h$ over one $\alpha$-fixed point is trivial. In particular, for $h$ odd, we may pick the $\alpha$-fixed point $$E=\left(\mathcal O\oplus \alpha \ldots \oplus \alpha^{\delta-1}\right)^{\oplus hr/\delta}\in SU_X(hr),$$ where $\delta$ stands for the order of $\alpha$.  The restriction of the tensor product map $t$ to $\{E\}\times A$ becomes $$t_{\alpha}:A\to U_X(hr), \,\,L\mapsto \left(\bigoplus_{i=0}^{\delta-1} L\otimes \alpha^i\right)^{\oplus hr/\delta}.$$ It suffices to verify that \begin{equation}\label{eqv1}t_{\alpha}^{\star}\Theta^h\cong h^{\star} \Theta^{r}\end{equation} holds $\alpha$-equivariantly. Note that the left hand side carries an $\alpha$-action since the morphism $t_{\alpha}$ factors through the quotient $$p:A\to A/\langle \alpha \rangle.$$ 

Let $$\pi:Y\to X$$ be the \'etale cover determined by $\alpha$; the curve $Y$ will be described in more detail in Section \ref{tracecalc} below. The kernel of the pullback $$\pi^{\star}:\text{Jac}(X)\to \text {Jac}(Y)$$ is generated by $\alpha$. Thus, we can write $$\pi^{\star}=i\circ p,$$ where $$i:\text {Jac}(X)/\langle \alpha \rangle\to \text {Jac} (Y)$$ denotes the inclusion. Writing $$\eta:\text{Jac}(Y)\to \text{Jac }(X)$$ for the norm map, we have $$\eta\circ \pi^{\star}=\delta, \text { hence } h/\delta\circ \eta\circ i \circ p= h.$$ 

Consider the morphisms $$g:U_X(\delta)\to U_X(hr), M\mapsto M^{\oplus hr/\delta}$$ and $$T_{\alpha}:\text{Jac}(Y)\to U_X(\delta),\,\,T_{\alpha}(M)=\pi_{\star} M.$$ Along the image of $\pi^{\star}$, we have $$T_{\alpha}(\pi^{\star}L)=\pi_{\star}\pi^{\star} L=\bigoplus_{i=0}^{\delta-1} L\otimes \alpha^i,$$ for $L\in \text{Jac} (X)$. Therefore, we conclude that $$t_{\alpha}=g\circ T_{\alpha}\circ i\circ p.$$ 

To establish that \eqref{eqv1} holds $\alpha$-equivariantly it suffices to confirm that $$i^{\star}T_{\alpha}^{\star}g^{\star} \Theta^h=i^{\star}\eta^{\star}\left(h/\delta\right)^{\star}\Theta^r.$$ In turn, this follows from the equalities $$g^{\star}\Theta_{\kappa}=\Theta^{hr/\delta}_{\kappa},\,\, T_{\alpha}^{\star}\Theta_{\kappa}=\Theta_{\pi^{\star}\kappa},\,\,$$ $$ i^{\star}\eta^{\star}\Theta_{\kappa}=i^{\star}\Theta^{\delta}_{\pi^{\star}\kappa}.$$ The last isomorphism is well known \cite {LB}, but for completeness let us briefly justify it here. First, the pullbacks of both sides via $p^{\star}$ are seen to agree. Indeed, we compute $$p^{\star} i^{\star} \eta^{\star} \Theta_{\kappa}=\delta^{\star} \Theta_{\kappa}=\Theta_{\kappa}^{\delta^2}\text{ and } $$ $$p^{\star} i^{\star}\Theta_{\pi^{\star}\kappa}^{\delta}=(\pi^{\star})^{\star}\Theta_{\pi^{\star}\kappa}^{\delta}=\Theta_{\pi_{\star}\pi^{\star}\kappa}^{\delta}=\otimes_{i=0}^{\delta-1}\Theta_{\kappa\otimes \alpha^i}^{\delta}=\Theta^{\delta^2}.$$ Thus, we conclude that $$i^{\star}\eta^{\star}\Theta_{\kappa}=i^{\star}\Theta^{\delta}_{\pi^{\star}\kappa}\otimes S$$ for some line bundle $S$ on $\text {Jac}(X)/\langle \alpha\rangle$ such that $p^{\star}S=\mathcal O.$ In particular, since $p$ has degree $\delta$, we have $S^{\delta}=\mathcal O.$ Now, it remains to observe that both $i^{\star}\eta^{\star}\Theta_{\kappa}$ and $i^{\star}\Theta^{\delta}_{\pi^{\star}\kappa}$ are both invariant under the action of the morphism $-1$, and thus $S$ must be invariant as well. In particular, since $\delta$ is odd, $S$ must be trivial.
The proof is now complete. 

\end{remark}

\section {The trace calculation} \label{tracecalc} The current section is central to the paper and is devoted to the proof of Theorem \ref{prop1}. We follow the strategy of \cite {B} and \cite {NR}. As a first step, we use the Hecke correspondence to transfer the calculation on a smooth moduli space. There, we determine the trace by means of Lefschetz-Riemann-Roch. Note that due to the involved geometry of the fixed loci, the details of the computation are significantly different than in the cases considered in \cite {B}, \cite {NR}.  
  
\subsection{Hecke correspondences} \subsubsection{} \label{hecke1}Theorem \ref{prop1} expresses the trace of a $\delta$-torsion line bundle ${\alpha}$ on the curve $X$ acting on the space of generalized theta functions $$H^{0}\left(SU_{X}(hr), \mathcal L^{hk}\right).$$ Here, the action of $\alpha$ on $SU_X(hr)$ is given by tensoring, while the line bundle $\mathcal L^{hk}$ is equipped with the $A[h]$-linearization of Remark \ref{rem1}. 

The trace of $\alpha$ is given in terms of Verlinde numbers of smaller rank on a higher genus curve $Y$. The curve $Y$ was already constructed in the previous section as a degree $\delta$ \'etale cover $$\pi:Y\to X$$ determined by the torsion point $\alpha.$ In fact, $$Y=\{s\in \alpha: s^{h}=1\}.$$ Since $\alpha$ has order precisely $\delta$, $Y$ is connected. Moreover, $$\text {genus of } Y=(g-1)\delta+1.$$ The Galois group $G$ of the cover is cyclic of order $\delta$. The line bundle $\alpha$ corresponds to a generator of the group $\widehat G$ of characters of $G$. Explicitly, we have $$\alpha=Y\times _{G} \mathbb C$$ where $G$ acts on the second factor with character $\alpha$. 
 
\subsubsection{} We transfer the trace computation on a smooth moduli space. Let $p\in X$, and let $SU_{X}(hr,p)$ be the moduli space of rank $hr$ bundles with determinant $\mathcal O_{X}(p)$. We have the following Hecke diagram \begin{center} $\xymatrix{  & \mathcal H \ar[dl]_{q} \ar[rd]^{q'} & \\ SU_{X}(hr) & &SU_{X}(hr,p),}$ \end{center} where $\mathcal H$ parametrizes pairs $$(F,E)\in SU_{X}(hr)\times SU_{X}(hr, p)\text{ such that } F\subset E.$$ 

As explained in \cite {BS}, pulling back and pushing forward by means of the above diagram induces an isomorphism \begin{equation}\label{isom1}H^{0}(SU_{X}(hr), \mathcal L^{hk})\cong H^{0}(SU_{X}(hr, p), \text {Sym}^{hk} \mathcal E_{p}).\end{equation} Here $$\mathcal E \to X\times SU_{X}(hr,p)$$ is the universal bundle, suitably normalized, and $\mathcal E_{p}$ is the restriction to $\{p\}\times SU_{X}(hr,p)$. The normalization of $\mathcal E$ is chosen such that \begin{equation}\label{normep}\det \mathcal E_{p}=\Theta_{A},\end{equation} where $\Theta_{A}$ is the ample generator of the Picard group of $SU_{X}(hr,p)$. The reference vector bundle $A$ on $X$ has rank $hr$ and slope $$\mu(A)=g-1-\frac{1}{hr}.$$ Note that the natural projection $q':\mathcal H\to SU_X(hr, p)$ has the structure of a projective bundle $\mathcal H\cong \mathbb P(\mathcal E_p)$ such that for the chosen normalization \eqref{normep} we have $$q^{\star} \mathcal L\cong \mathcal O_{\mathcal H}(1).$$ Pushing forward by $q'$ we obtain $$q'_{\star}q^{\star}\mathcal L^{hk}\cong \text{Sym}^{hk} \mathcal E_p,$$ which implies \eqref{isom1}.

The isomorphism \eqref{isom1} is invariant under the action of $\alpha$ \cite {B}. In Subsection \ref{hecke1} we explained how $\alpha$ acts on the left hand side. As for the right hand side, the action of $\alpha\in \widehat G$ on the moduli space $SU_{X}(hr,p)$ is given, as usual, by tensoring. The action of $\alpha$ on $\text{Sym}^{hk}\mathcal E_{p}$ is obtained as follows. We have isomorphisms $$\alpha^{\star}\mathcal E \cong\mathcal E\otimes \alpha$$ unique up to a scalar \cite{NR}, inducing isomorphisms $$u:\alpha^{\star}\mathcal E_{p}\to \mathcal E_{p}.$$ The requirement $u^{h}=1$ fixes $u$ up to an $h$-root of unity. The induced isomorphism $\text {Sym}^{hk}u$ on $\text{Sym}^{hk}\mathcal E_{p}$ is unambiguously defined. The identification $$q^{\star}\mathcal L^h\cong \mathcal O_{\mathcal H}(h)$$ is $\widehat G$-equivariant since the $\widehat G$-action on both sides is trivial in the fibers over the $\alpha$-fixed loci. Therefore, we have \begin{equation}\label{isom2}\text{Trace}\left({\alpha}, H^{0}\left(SU_{X}(hr), \mathcal L^{hk}\right)\right)=\text{Trace}\left({\alpha}, H^{0}\left(SU_{X}(hr, p), \text {Sym}^{hk}\mathcal E_{p}\right)\right).\end{equation}

\subsubsection{} Similarly, let $\mathcal M$ be the ample generator of the Picard group of $SU_{Y}\left(\frac{hr}{\delta}\right)$. The same argument shows that \begin{equation}\label{isom3}H^{0}\left(SU_{Y}\left(\frac{hr}{\delta}\right), \mathcal M^{\frac{hk}{\delta}}\right)=H^{0}\left(SU_{Y}\left(\frac{hr}{\delta}, q\right), \text{Sym}^{\frac{hk}{\delta}}\,\widetilde {\mathcal D}_{q}\right),\end{equation} where $q\in Y$ and $$\widetilde {\mathcal D}\to SU_{Y}\left(\frac{hr}{\delta}, q\right)\times Y$$ is the universal bundle, suitably normalized. In fact, we require \begin{equation}\label{normalizeds}\det \widetilde {\mathcal D}_{q}=\Theta_{B},\end{equation} where $\Theta_{B}$ stands for the ample generator of the Picard group of $SU_{Y}\left(\frac{hr}{\delta}, q\right)$. The corresponding reference bundle $B$ on $Y$ has rank $\frac{hr}{\delta}$ and slope $$\mu(B)=\delta(g-1)-\frac{\delta}{hr}.$$

Combining \eqref{isom2} and \eqref{isom3}, we see that the Theorem is equivalent to the equality \begin{equation}\label{reform}\text{Trace}\left({\alpha}, H^{0}\left(SU_{X}(hr, p), \text {Sym}^{hk}\mathcal E_{p}\right)\right)=\left({\frac{r}{r+k}}\right)^{-(g-1)(\delta-1)}\end{equation} $$\times h^{0}\left(SU_{Y}\left(\frac{hr}{\delta}, q\right), \text{Sym}^{\frac{hk}{\delta}}\,\widetilde {\mathcal D}_{q}\right).$$

\subsection {Lefschetz-Riemann-Roch} We will verify \eqref{reform} using Lefschetz-Riemann-Roch.

\subsubsection{} We first describe the fixed locus of the action of $\alpha$ on $SU_{X}(hr,p)$. Pushforward by $\pi:Y\to X$ induces a morphism $$\pi_{\star}:U_{Y}\left(\frac{hr}{\delta}, 1\right)\to U_{X}(hr, 1), E\mapsto \pi_{\star}E.$$ Note that stability of bundles is preserved under $\pi_{\star}$, as shown in \cite {NR}. Therefore, we obtain a commutative diagram \begin{center} $\xymatrix{  U_{Y}\left(\frac{hr}{\delta}, 1\right) \ar[r]^{\pi_{\star}} \ar[d]^{\det} & U_{X}(hr,1) \ar [d]^{\det} \\ \text{Jac}^{1}(Y)\ar[r]^{\eta} &\text{Jac}^{1}(X),}$ \end{center} where $\eta$ is the norm map. It is shown in \cite {NR} that the $\alpha$-fixed locus on $U_{X}(hr,1)$ coincides with the image of the morphism $\pi_{\star}$. Restricting to $SU_{X}(hr, p),$ the fixed locus is the subscheme $$i:SU_{Y}\left(\frac{hr}{\delta}, P\right)\hookrightarrow SU_{X}(hr,p)$$ consisting of those bundles $E$ with the requirement that the determinants $$\det E \in P.$$ Here $P$ is a connected component of $\eta^{-1}(p)\hookrightarrow \text {Jac}^{1}(Y).$ We clearly have $$\dim {P}=(g-1)(\delta-1).$$

\subsubsection{} By Lefschetz-Riemann-Roch, we obtain \begin{equation}\label{lrr}\text{Trace}\left({\alpha}, H^{0}\left(SU_{X}(hr, p),\text{Sym}^{hk}{\mathcal E_{p}}\right)\right)\end{equation} $$=\int_{SU_{Y}\left(\frac{hr}{\delta}, P\right)}i^{\star}{\text {ch}}_{\widehat G}\left(\text{Sym}^{hk}{\mathcal E_{p}}\right)(\alpha)\,\,\left(\prod_{x\neq 1} \text{ch}_{-\langle\alpha^{-1},x\rangle}N_{x}^{\vee}\right)^{-1} \text{Todd}\left(SU_{Y}\left(\frac{hr}{\delta}, P\right)\right).$$ To explain the notation, the first Chern character is considered $\widehat G$-equivariantly as an element in $H^{\star}_{\widehat G}\left(SU_{Y}\left(\frac{hr}{\delta}, P\right)\right)\cong H^{\star}\left(SU_{Y}\left(\frac{hr}{\delta}, P\right)\right)\otimes R(\widehat G)$, and then evaluated against the element $\alpha\in \widehat G$. Also, we write $N$ for the normal bundle of the fixed locus $$i:SU_{Y}\left(\frac{hr}{\delta}, P\right)\hookrightarrow SU_{X}(hr,p).$$ For $x \in G$, $N_{x}$ denotes the $x$-eigensubbundle of $N$ corresponding to the character of $\widehat G$ determined by $x\in G$. Finally, we set $$\text{ch}_{t}(V)=\prod(1+te^{x_{i}})$$ for a bundle $V$ with Chern roots $x_{i}$. 

\subsection{The intersection computation} We proceed to compute the intersection number \eqref{lrr}. We evaluate the three terms appearing in the integral \eqref{lrr} one by one. Then, we use an \'etale pullback to split off an abelian factor, thereby bringing the intersection number \eqref{lrr} in the form predicted by \eqref{reform}. 

\subsubsection{The Chern character} We evaluate first the equivariant Chern character appearing in \eqref{lrr}. Let $$\mathcal D\to SU_{Y}\left(\frac{hr}{\delta}, P\right)\times Y$$ be an universal bundle. We let $$\pi^{-1}(p)=\{p_{1}, \ldots, p_{\delta}\}$$ be the preimages of $p$, which differ by the action of $G$. Possibly after relabeling, we may assume that $p_1\in P$. We let $\mathcal D_{p_{i}}$ be the restriction of $\mathcal D$ to $\{p_{i}\}\times SU_{Y}\left(\frac{hr}{\delta}, P\right)$. It is shown in \cite {B} that after a suitable normalization of $\mathcal D$ we may assume that \begin{equation}\label{ed}\mathcal E_{p}\cong \bigoplus_{i=1}^{\delta}\mathcal D_{p_{i}}.\end{equation}  We will check later that this normalization of $\mathcal D$ agrees with \eqref{normalizeds}. Since the bundles $\mathcal D_{p_{i}}$ are algebraically equivalent, they have the same Chern classes. Therefore, for the chosen normalization of $\mathcal D$ we have \begin{equation}\label{firstc}c_{1}(\mathcal D_{p_{1}})=\frac{1}{\delta}i^{\star}c_{1}(\mathcal E_{p})=\frac{1}{\delta}i^{\star}c_{1}(\Theta_{A}).\end{equation}

As shown in \cite {B} or section $10.4$ in \cite {BLS}, equation \eqref{ed} gives the $\widehat G$-eigenbundle decomposition of $\mathcal E_{p}$. As a consequence, $${\text {ch}}_{\widehat G}\left(\text{Sym}^{hk}{\mathcal E_{p}}\right)(\alpha)=\text{ch}_{\widehat G}\left(\text{Sym}^{hk}{(\oplus_{i=1}^{\delta} \mathcal D_{p_{i}})}\right)(\alpha).$$
The corresponding eigenbundles are in bijection with the elements $g\in G$. The Chern character evaluated against $\alpha$ is an expression in the $\delta$-roots of unity $\langle \alpha, g\rangle.$ We claim that this expression does not depend on the choice of $\alpha$.  
Indeed, letting $$\theta_{1}, \ldots, \theta_{\frac{hr}{\delta}}$$ denote the Chern roots of the bundle $\mathcal D_{p_{1}},$ the $\widehat {G}$-equivariant Chern roots of the bundle $\oplus_{i=1}^{\delta}\mathcal D_{p_{i}}$ are of the form $$\theta_{1}\otimes g, \ldots, \theta_{\frac{hr}{\delta}}\otimes g$$ for $g\in G$. We form the generating series $$\sum_{m}\text{ch}_{\widehat G}\left(\text{Sym}^{m}({\oplus_{i=1}^{\delta} \mathcal D_{p_{i}}})\right)t^{m}=\prod_{j=1}^{\frac{hr}{\delta}}\prod_{g\in G}\frac{1}{1- (e^{\theta_{j}}\otimes g)\,t}.$$ Evaluating against $\alpha$ gives $$\sum_{m}\text{ch}_{\widehat G}\left(\text{Sym}^{m}(\oplus_{i=1}^{\delta} \mathcal D_{p_{i}})\right)(\alpha) \,t^{m}=\prod_{j=1}^{\frac{hr}{\delta}}\prod_{g\in G} \frac{1}{1-\langle \alpha, g\rangle e^{\theta_{j}}t}=\prod_{j=1}^{\frac{hr}{\delta}}\frac{1}{1-e^{\delta \theta_{j}}\,t^{\delta}}.$$  This shows that \begin{equation}\label{fund1}\text{ch}_{\widehat G}\left(\text{Sym}^{hk}{\mathcal E_{p}}\right)(\alpha)=\text{ch}\left(\text{Sym}^{\frac{hk}{\delta}}{\mathcal F}\right),\end{equation} where $\mathcal F$ is a bundle on $SU_Y(\frac{hr}{\delta}, P)$ with Chern roots $\delta \theta_{1}, \ldots, \delta \theta_{\frac{hr}{\delta}}.$ 

\subsubsection{The equivariant normal bundles} Next, we compute the expression $$\prod_{x\in G\setminus\{1\}}\text {ch}_{-\langle \mu, x\rangle} N_{x}^{\vee}$$ for a generator $\mu\in \widehat G$. The bundles $N_x$ are described in \cite {NR}. We have $$N_{x}=R^{1}p_{!}(V_{x})$$ where $V_{x}$ is a vector bundle on $SU_{Y}\left(\frac{hr}{\delta}, P\right)\times X,$ and $$p:SU_{Y}\left(\frac{hr}{\delta}, P\right) \times X\to SU_{Y}\left(\frac{hr}{\delta},P\right)$$ is the projection. The pullback of $V_{x}$ via $$\pi:SU_{Y}\left(\frac{hr}{\delta}, P\right)\times Y\to SU_{Y}\left(\frac{hr}{\delta}, P\right)\times X$$ is given by \begin{equation}\label{pullback}\pi^{\star}V_{x}=\bigoplus_{g\in G}g^{\star}\left(x^{\star}\mathcal D\otimes \mathcal D^{\vee}\right).\end{equation} Similarly, we have $$\pi^{\star}V_{x^{-1}}=\bigoplus_{g\in G}g^{\star}\left((x^{-1})^{\star}\mathcal D\otimes \mathcal D^{\vee}\right)=\bigoplus_{g\in G}g^{\star}x^{\star}\left((x^{-1})^{\star}\mathcal D\otimes \mathcal D^{\vee}\right)$$ $$=\bigoplus_{g\in G}g^{\star}\left(\mathcal D\otimes x^{\star}\mathcal D^{\vee}\right)=\pi^{\star}V_{x}^{\vee}.$$ Since $\pi$ is injective in cohomology, it follows that \begin{equation}\label{eqal}\text{ch}(V_{x})=\text{ch}(V_{x^{-1}}^{\vee}).\end{equation} We restrict \eqref{pullback} to $SU_Y\left(\frac{hr}{\delta}, P\right)\times \{\text{pt}\}$, recalling the remark preceding \eqref {firstc}. We obtain \begin{equation}\label{eqal1}\text {ch }V_{x}|_{SU_{Y}\times\{\text{pt}\}}=\delta \text { ch }(\mathcal D\otimes \mathcal D^{\vee})|_{SU_{Y}\times \{\text{pt}\}}.\end{equation} Here, we used the obvious shorthand $SU_Y$ for the moduli space $SU_Y\left(\frac{hr}{\delta}, P\right)$. 

Grothendieck-Riemann-Roch for the projection $p$ gives \begin{equation}\label{w1}\text {ch } N_{x}=-p_{!}\left(\left(1-(g-1)\omega\right)\text{ch } V_{x}\right)=(g-1)\text{ch }V_{x}{{|}}_{SU_{Y}\times\{\text{pt}\}}-p_{!}\text{ ch } V_{x},\end{equation} where $\omega$ is the class of a point on $X$. Similarly, \begin{equation}\label{w2}\text{ch }N_{x^{-1}}=(g-1)\text{ch }V_{x^{-1}}{{|}}_{SU_{Y}\times\{\text{pt}\}}-p_{!}\text { ch }V_{x^{-1}}.\end{equation} Thus, using \eqref{eqal}, \eqref{eqal1}, \eqref{w1} and \eqref{w2} we obtain \begin{equation}\label{chernc}\text{ch }N_{x}+\text{ch }N_{x^{-1}}^{\vee}=2(g-1)\delta \text { ch }(\mathcal D\otimes \mathcal D^{\vee})|_{SU_{Y}\times\{\text{pt}\}}.\end{equation} Equation \eqref{chernc} needs to be checked separately for the even and odd pieces of the Chern character. 

Let us partition the elements of $G\setminus\{1\}$ into two sets $$G\setminus\{1\}=\{x_{1}, \ldots, x_{l}\}\cup \{x_{1}^{-1}, \ldots, x_{l}^{-1}\}$$ such that $$\langle \mu, x_{i}^{-1}\rangle=\exp\left(\frac{2\pi i\sqrt{-1}}{\delta}\right) \text{for } 1\leq i\leq l,$$ where $\delta=2l+1$. This partitioning depends on the choice of $\mu$. We have $$\prod_{x\neq 1}\text {ch}_{-\langle \mu, x\rangle}N_{x}^{\vee}=\prod_{i=1}^{l} \text {ch}_{-\langle \mu, x_{i}\rangle} N^{\vee}_{x_{i}}\cdot \text{ch}_{-\langle \mu, x_{i}^{-1}\rangle} N^{\vee}_{x_{i}^{-1}}.$$ Now, for any vector bundle $V$, we have the identity \begin{equation}\label{chernc1}\text{ch}_{-t} V^{\vee}=(-t)^{\text {rk}V}e^{-c_{1}(V)}\text {ch}_{-t^{-1}} V.\end{equation} Therefore, by \eqref{chernc} and \eqref{chernc1},
$$\prod_{x\neq 1}\text {ch}_{-\langle \mu, x\rangle}N_{x}^{\vee}= \prod_{i=1}^{l}(-\langle\mu,x_{i}\rangle)^{\frac{h^{2}r^2}{\delta}(g-1)} e^{-c_{1}(N_{x_{i}})} \text {ch}_{-\langle \mu, x_{i}^{-1}\rangle} N_{x_{i}}\cdot \text{ch}_{-\langle \mu, x_{i}^{-1}\rangle} N_{x_{i}^{-1}}^{\vee}$$ $$=(-1)^{lr(g-1)}e^{-c_{1}(N_{x_{1}})-\ldots-c_{1}(N_{x_{l}})} \prod_{i=1}^{l} \text {ch}_{-\exp\left(\frac{2\pi i\sqrt{-1}}{\delta}\right)}\left(\mathcal D\otimes \mathcal D^{\vee}|_{SU_{Y}\times \{{\text{pt}\}}}\right)^{2(g-1)\delta}.$$ In the above, we made use of the fact that $\delta|h$ and that $h$ is odd to pin down the factor $$  (-\langle\mu,x_{i}\rangle)^{\frac{h^{2}}{\delta}(g-1)}=(-1)^{g-1}.$$ Using \eqref{chernc1} again, we see that the expression $$\prod_{i=1}^{l}\text {ch}_{-\exp\left(\frac{2\pi i\sqrt{-1}}{\delta}\right)}\left(\mathcal D\otimes \mathcal D^{\vee}|_{SU_{Y}\times \{{\text{pt}\}}}\right)^{2(g-1)\delta}$$ equals $$(-1)^{lr(g-1)}\prod_{i=1}^{l}\left\{\text{ch}_{-\exp\left(\frac{2\pi i\sqrt{-1}}{\delta}\right)}\left(\mathcal D\otimes \mathcal D^{\vee}|_{SU_{Y}\times \{{\text{pt}\}}}\right) \text{ch}_{-\exp\left(-\frac{2\pi i\sqrt{-1}}{\delta}\right)}\left(\mathcal D\otimes \mathcal D^{\vee}|_{SU_{Y}\times \{{\text{pt}\}}}\right)\right\}^{(g-1)\delta}$$ \begin{equation}\label{prodf}=(-1)^{lr(g-1)}\prod_{i=1}^{\delta-1}\left\{\text{ch}_{-\exp\left(\frac{2\pi i\sqrt{-1}}{\delta}\right)}\left(\mathcal D\otimes \mathcal D^{\vee}|_{SU_{Y}\times \{{\text{pt}\}}}\right)\right\}^{(g-1)\delta}.\end{equation} It is easy to see that for any line bundle $L$ we have $$\prod_{i=1}^{\delta-1} \text{ch}_{-\exp\left(\frac{2\pi i\sqrt{-1}}{\delta}\right)}(L)=\frac{1-\exp(\delta c_1(L))}{1-\exp (c_1(L))}.$$ Now, recall that we wrote $$\theta_{1}, \ldots, \theta_{\frac{hr}{\delta}}$$ for the Chern roots of $\mathcal D|_{SU_{Y}\times \{\text{pt}\}}.$ Then expression \eqref{prodf} becomes $$(-1)^{lr(g-1)}\prod_{i,j}\left(\frac{1-\exp(\delta\theta_{i}-\delta\theta_{j})}{1-\exp(\theta_{i}-\theta_{j})}\right)^{(g-1)\delta}=(-1)^{lr(g-1)}\delta^{hr(g-1)}\prod_{i<j}\left(\frac{\sinh\frac{\delta\theta_{i}-\delta\theta_{j}}{2}}{\sinh\frac{\theta_{i}-\theta_{j}}{2}}\right)^{2(g-1)\delta}.$$ The prefactor $\delta^{hr(g-1)}$ corresponds to the case $i=j$. 

Now, from \eqref{chernc} we have $$c_{1}(N_{x})=c_{1}(N_{x^{-1}}).$$ Therefore, $$c_{1}(N_{x_{1}})+\ldots+c_{1}(N_{x_{l}})=\frac{1}{2}\sum_{x\neq 1}c_{1}(N_{x})=\frac{1}{2}\left(i^{\star}c_{1}(SU_{X}(hr,p))-c_{1}\left(SU_{Y}\left(hr/\delta, P\right)\right)\right).$$ Moreover, from \cite {DN}, we know that $$\frac{1}{2}c_{1}(SU_{X}(hr,p))=\Theta_{A}.$$

Putting everything together, we conclude that \begin{equation}\label{fund2}\prod_{x\neq 1}\text {ch}_{-\langle \mu, x\rangle}N_{x}^{\vee}=\delta^{hr(g-1)}\exp\left(-i^{\star}\Theta_{A}+\frac{1}{2}c_{1}\left(SU_{Y}\left(hr/\delta, P\right)\right)\right)\prod_{i<j}\left(\frac{\sinh\frac{\delta\theta_{i}-\delta\theta_{j}}{2}}{\sinh\frac{\theta_{i}-\theta_{j}}{2}}\right)^{2(g-1)\delta}.\end{equation} 

\subsubsection {The Todd character} Finally, it is proved in \cite {N} that \begin{equation}\label{fund3}\text {Todd}\left(SU_{Y}\left(hr/\delta, P\right)\right)=\exp\left(\frac{1}{2}c_{1}\left(SU_{Y}\left({hr}/{\delta}, P\right)\right)\right) \prod_{i<j}\left(\frac{\theta_{i}-\theta_{j}}{2\sinh \frac{\theta_{i}-\theta_{j}}{2}}\right)^{2(g-1)\delta}.\end{equation} Collecting all terms appearing in \eqref{lrr} with the aid of \eqref{fund1}, \eqref{fund2} and \eqref{fund3}, we obtain that the trace of $\alpha$ is \begin{equation}\label{beforepb}\delta^{-\frac{h^{2}r^2}{\delta}(g-1)}\int_{SU_{Y}\left(\frac{hr}{\delta}, P\right)}\exp\left(i^{\star}c_{1}(\Theta_{A})\right)\text{ch}\left(\text{Sym}^{\frac{hk}{\delta}}{\mathcal F}\right) \prod_{i<j}\left(\frac{\delta\theta_{i}-\delta\theta_{j}}{2\sinh \frac{\delta \theta_{i}-\delta\theta_{j}}{2}}\right)^{2(g-1)\delta}.\end{equation}

\subsubsection {Splitting off the abelian factor} Next, we will fix the determinant of the bundles in the moduli space on which \eqref{beforepb} is evaluated. To this end, note that twisting by the point $p_1\in P$ gives an isomorphism $$j:\text{Jac}^{1}(Y)\to \text{Jac}(Y).$$ The image $Q=j(P)$ becomes an abelian subvariety of $\text {Jac}(Y)$. We will make use of the following fiber diagram \begin{center}$\xymatrix{SU_{Y}\left(\frac{hr}{\delta}, p_{1}\right)\times Q \ar[r]^{t} \ar[d]^{p} & SU_{Y}\left(\frac{hr}{\delta}, P\right) \ar[d]^{j\circ\det}\ar[r]^{i}& SU_{X}(hr, p) \\ Q\ar[r]^{\frac{hr}{\delta}} &  Q &}$ \end{center} where $t$ is the tensor product map.  

We will compute the pullbacks under $t$ of various universal structures. We claim that \begin{equation}\label{restr}t^{\star}i^{\star}\Theta_{A}\cong\Theta_B^{\delta}\boxtimes \Theta^{\frac{h^{2}r^2}{\delta}},\end{equation} where the second factor $\Theta$ is the restriction of a Theta bundle to $Q\hookrightarrow \text{Jac}(Y)$. To see the last isomorphism, observe first that $t^{\star}i^{\star}\Theta_{A}$ is supported on the locus $$\{(E,L): h^{0}(E\otimes L\otimes \pi^{\star}A)\neq 0\}.$$ Its restriction to $SU_Y\left(\frac{hr}{\delta}, p_{1}\right)\times \{L\}$ is independent of $L\in Q$; in fact it equals $\Theta_{L\otimes \pi^{\star}A}=\Theta_B^{\delta}$ by the results of \cite {DN}. Equation \eqref{restr} follows from the see-saw theorem.

Let us write $\widetilde {\mathcal D}$ for the restriction of $\mathcal D$ to $$SU_{Y}\left(\frac{hr}{\delta}, p_{1}\right)\times Y\hookrightarrow SU_{Y}\left(\frac{hr}{\delta}, P\right)\times Y.$$ The bundle $\widetilde {\mathcal D}$ serves as a universal bundle on $SU_{Y}\left(\frac{hr}{\delta}, p_{1}\right)\times Y$. We claim that on $SU_Y\left(\frac{hr}{\delta}, p_{1}\right)\times Q$ we have \begin{equation}\label{pullt}t^{\star}\mathcal D_{p_{1}}\cong \widetilde{\mathcal D}_{p_{1}}\boxtimes L,\end{equation} for some line bundle $L$ on $Q$. To check this isomorphism, first observe that the vector bundles $$(t\times 1)^{\star}\mathcal D \text { and } \widetilde{\mathcal D}\otimes \mathcal P$$ on $SU_{Y}\left(\frac{hr}{\delta}, p_{1}\right)\times Q \times Y$ have the same restriction to $\{E\}\times \{q\}\times Y,$ for $E\in SU_{Y}\left(\frac{hr}{\delta}, p_{1}\right)$, $q\in Q$. Here, $\mathcal P$ is the Poincare bundle on $Q\times Y\hookrightarrow \text {Jac }(Y)\times Y$. Therefore, by lemma $2.5$ of \cite {R}, we can write $$(t\times 1)^{\star}\mathcal D\cong\widetilde{\mathcal D}\otimes \mathcal P\otimes \widetilde {L}$$ for a line bundle $\widetilde L$ on $SU_{Y}\left(\frac{hr}{\delta}, p_{1}\right)\times Q$. Restricting to $SU_{Y}\left(\frac{hr}{\delta}, p_{1}\right)\times Q\times \{p_{1}\}$, we obtain \begin{equation}\label{pullt1}t^{\star}\mathcal D_{p_{1}}\cong \widetilde{\mathcal D}_{p_{1}}\otimes \widetilde L',\end{equation} for some line bundle $\widetilde L'$ on $SU_{Y}\left(\frac{hr}{\delta}, p_{1}\right)\times Q$. To prove \eqref{pullt}, it suffices to show that $\widetilde L'$ is the pullback of a line bundle $L$ on $Q$. In turn, we will show that the restriction of $\widetilde L'$ to $SU_{Y}\left(\frac{hr}{\delta}, p_{1}\right)\times \{q\}$ is trivial, for all $q\in Q$. Since $SU_{Y}\left(\frac{hr}{\delta}, p_{1}\right)$ is simply connected, it remains to verify that $$c_{1}(\widetilde L'|_{SU_{Y}\times \{q\}})=0.$$ Let $t_{q}$ denote the restriction of $t$ to $SU_{Y}\left(\frac{hr}{\delta}, p_{1}\right)\times \{q\}$.  From \eqref{pullt1}, we compute 
 $$\frac{hr}{\delta}c_{1}(\widetilde L'|_{SU_{Y}\times \{q\}})=c_{1}(t_{q}^{\star}\mathcal D_{p_{1}})-c_{1}(\widetilde{\mathcal D}_{p_{1}}).$$ The last expression vanishes since by continuity $c_{1}(t_{q}^{\star}D_{p_{1}})$ is independent of $q\in Q$. This proves \eqref{pullt}.
  
We will identify the line bundle $L$. Taking determinants in \eqref{pullt} and using \eqref{firstc}, we have $$c_{1}(\widetilde{\mathcal D}_{p_{1}})+\frac{hr}{\delta}c_{1}(L)=t^{\star}c_{1}(\mathcal D_{p_{1}})=\frac{1}{\delta}t^{\star}i^{\star}c_{1}(\Theta_{A}).$$ Comparing with \eqref{restr}, we have $$c_{1}(\widetilde{\mathcal D}_{p_{1}})=\Theta_{B},\,\, c_{1}(L)=\frac{hr}{\delta}c_{1}(\Theta).$$ In particular, since $SU_{Y}\left(\frac{h}{\delta}, p_{1}\right)$ is simply connected, we have \begin{equation}\label{norm}\det \widetilde{\mathcal D}_{p_{1}}\cong \Theta_{B}.\end{equation} This confirms the normalization \eqref{normalizeds}. Since we work cohomologically, somewhat abusively we may replace $L$ by $\Theta^{\frac{hr}{\delta}}.$ Therefore, from \eqref{pullt} we have $$t^{\star}\mathcal F\cong \widetilde {\mathcal F}\boxtimes \Theta^{hr}$$ where $\widetilde {\mathcal F}$ is a bundle with Chern roots $\delta\tilde\theta_{1}, \ldots, \delta\tilde\theta_{\frac{hr}{\delta}}$. Here, we wrote $\tilde \theta_{1}, \ldots, \tilde \theta_{\frac{hr}{\delta}}$ for the Chern roots of $\widetilde {\mathcal D}_{p_{1}}$. We conclude \begin{equation}\label{pullt2}t^{\star}\text {Sym}^{\frac{hk}{\delta}}\mathcal F=\text{Sym}^{\frac{hk}{\delta}}\widetilde {\mathcal F}\boxtimes \Theta^{\frac{h^{2}rk}{\delta}}.\end{equation}

Pulling back the intersection number \eqref{beforepb} under $t$ with the aid of \eqref{restr}, \eqref{pullt} and \eqref{pullt2}, we obtain that $\text{Trace } \alpha$ equals \begin{equation}\label{trace2} \delta^{-\frac{h^{2}r^2}{\delta}(g-1)}\int_{SU_{Y}\left(\frac{hr}{\delta}, p_{1}\right)}\text { ch }(\text {Sym}^{\frac{hk}{\delta}}\widetilde{\mathcal F})\, \exp(\delta \Theta_B)\, \prod_{i<j}\left(\frac{\delta\tilde\theta_{i}-\delta\tilde\theta_{j}}{2\sinh \frac{\delta \tilde\theta_{i}-\delta\tilde\theta_{j}}{2}}\right)^{2(g-1)\delta}\times \end{equation} $$\frac{1}{\deg t}\int_{Q}\exp\left(\frac{h^{2}}{\delta}r(r+k)\Theta\right).$$ 

\subsubsection {The trace evaluation} We will now evaluate the two integrals of \eqref{trace2}. The second integral equals $$\frac{1}{\deg t}\int_{Q}\exp\left(\frac{h^{2}}{\delta}r(r+k)\Theta\right)=\left(\frac{(r+k)\delta}{r}\right)^{\dim Q}\frac{1}{\deg (hr/\delta)}\int_{Q}\exp\left(\left(\frac{hr}{\delta}\right)^{\star}\Theta\right)$$  \begin{equation}\label{fund4}=\left(\frac{(r+k)\delta}{r}\right)^{(\delta-1)(g-1)}\chi(P, \Theta)= \left(\frac{(r+k)\delta}{r}\right)^{(\delta-1)(g-1)}\delta^{g-1}.\end{equation} The Euler characteristic $$\chi(P, \Theta)=\delta^{g-1}$$ is computed in Corollary $4.16$ of \cite {NR}. 

The first integral in \eqref{trace2} is calculated a backward application of Riemann-Roch. We have $$\int_{SU_{Y}\left(\frac{hr}{\delta}, p_{1}\right)}\text { ch }(\text {Sym}^{\frac{hk}{\delta}}\widetilde{\mathcal F}) \, \exp(\delta \Theta_B) \,\prod_{i<j}\left(\frac{\delta\tilde\theta_{i}-\delta\tilde\theta_{j}}{2\sinh \frac{\delta \tilde\theta_{i}-\delta\tilde\theta_{j}}{2}}\right)^{2(g-1)\delta}$$ $$=\delta^{\dim SU_{Y}\left(\frac{hr}{\delta},p_{1}\right)} \int_{SU_{Y}\left(\frac{hr}{\delta}, p_{1}\right)} \text {ch }(\text {Sym}^{\frac{hk}{\delta}} \,\widetilde{\mathcal D}_{p_{1}}) \,\exp\left(\frac{1}{2}c_{1}\left(SU_{Y}\right)\right) \prod_{i<j}\left(\frac{\tilde\theta_{i}-\tilde\theta_{j}}{2\sinh \frac{\tilde\theta_{i}-\tilde\theta_{j}}{2}}\right)^{2(g-1)\delta}.$$ Recalling \eqref{fund3}, this rewrites as
$$\delta^{(\frac{h^{2}r^2}{\delta^{2}}-1)\delta(g-1)} \int_{SU_{Y}\left(\frac{hr}{\delta}, p_{1}\right)}\text { ch }(\text {Sym}^{\frac{hk}{\delta}} \,\widetilde{\mathcal D}_{p_{1}})\text{ Todd } \left(SU_{Y}\left(\frac{hr}{\delta}, p_{1}\right)\right)$$
$$=\delta^{(\frac{h^{2}r^2}{\delta^{2}}-1)\delta(g-1)}\chi\left(SU_{Y}\left(\frac{hr}{\delta}, p_{1}\right), \text{Sym}^{\frac{hk}{\delta}}\, \widetilde{\mathcal D}_{p_{1}}\right).$$ Finally, substituting into \eqref{trace2} we obtain expression \eqref{reform}. This completes the proof of the Theorem.
\qed

\section {The Verlinde formula for $PGL_{r}$}

We will show how the trace computation of Theorem \ref{prop1} leads to the Verlinde formula for non-simply connected groups covered by $SL_{r}$. Similar results for the group $PGL_{r}$ with $r$ prime are found in \cite {B}, and for arbitrary non-simply connected groups in \cite{AMW}.

To fix notation, write $$\bfz_{r}\hookrightarrow SL_{r}$$ for the center of $SL_{r}$. For each $d|r$, we let $\bfz_{d}$ denote the subgroup of the center consisting of elements of order dividing $d$. Set $$\bfg_{d}=SL_{r}/\bfz_{d}.$$ The moduli space of $\bfg_{d}$-bundles has $d$ connected components indexed by elements in $\pi_{1}(\bfg_{d})\cong \mathbb Z/d\mathbb Z.$ As explained in \cite {BLS}, the irreducible component which contains the trivial bundle is the quotient $$\bfM_{d}=SU_{X}(r)/A[d].$$ As usual, the action of $A[d]$ on $SU_{X}(r)$ is obtained by tensoring.

Let $\mathcal G$ be the generator of the Picard group on the moduli stack of $\bfg_{d}$-bundles; note that $\mathcal G$ does not descend to the moduli scheme. Nonetheless, for $r$ odd, the line bundle $$\mathcal L^{d}\to SU_{X}(r)$$ descends to the line bundle $\mathcal G^{d}$ on the the moduli scheme $\bfM_{d}$, cf. Proposition $9.1$ in \cite {BLS}. Pick $k$ such that $d|k$. We have $$h^{0}(\bfM_{d}, \mathcal G^{k})=h^{0}(SU_{X}(r), \mathcal L^{k})^{A[d]}.$$ Thus, $$h^{0}(\bfM_{d}, \mathcal G^{k})=\frac{1}{d^{2g}}\sum_{\alpha\in A[d]} \text {Trace}\left({\alpha}, H^{0}(SU_{X}(r), \mathcal L^{k})\right)=$$ \begin{equation}\label{comba}=\frac{1}{d^{2g}} \frac{r^g}{(r+k)^{g}}\sum_{\delta|d}\bfn ({\delta})\,\bfv_{{\delta}(g-1)+1}\left(\frac{r}{\delta}, \frac{k}{\delta}\right).\end{equation} Here $\bfn(\delta)$ denotes the number of elements in the Jacobian of order precisely $\delta$. It follows from \eqref{sum} applied to the trivial character that $$\bfn(\delta)=\delta^{2g} \prod_{p|\delta} \left(1-\frac{1}{p^{2g}}\right).$$ We clearly must have \begin{equation}\label{sumbfn}\sum_{\delta|m}\bfn(\delta)=m^{2g},\end{equation} for any positive integer $m$. 

We will rewrite equation \eqref{comba} in a form that agrees with the more general Theorem $5.1$ in \cite {AMW}:

\begin{corollary}\label{combf}Assuming that $d$ divides $r$ and $k$, and $r$ is odd, we have $$h^{0}(\bfM_{d}, \mathcal G^{k})=r^{g}(r+k)^{(r-1)(g-1)-1}\sum_{\stackrel{S\subset\{1, \ldots, r+k\}}{|S|=r}} \xi_{d}(S) \prod_{s\neq t\in S}\left|2\sin \frac{s-t}{r+k} \pi \right|^{1-g}.$$ \end{corollary}

Here, we set $$\xi_{d}(S)=\left(\frac{\gcd(\delta_{S}, d)}{d}\right)^{2g},$$ where $\delta_{S}$ is the largest coperiod of the set $S$. Our terminology is as follows. A subset $S\subset \{1, \ldots, r+k\}$ with $r$ elements is called $\delta$-coperiodic if $$\delta|\gcd(r,k)$$ and there exists $$\widetilde S\subset \left\{1, \ldots, \frac{r+k}{\delta}\right\},\,\,|\widetilde S|=\frac{r}{\delta},$$ such that $$S=\widetilde S_{0}\cup \widetilde S_{1}\cup \ldots\cup \widetilde S_{\delta-1},$$ where $$\widetilde S_{i}=\frac{r+k}{\delta}\,i+\widetilde S, \,\, 0\leq i\leq \delta-1.$$ For example, any set $S\subset \{1, \ldots, r+k\}$ has coperiod $1$. 

Let $S$ be a $\delta$-coperiodic set and let $\widetilde S$ as above. We have $$\prod_{s\neq t\in S}\left | 2\sin \frac{s-t}{r+k}\pi \right |=\prod_{\tilde s, \tilde t\in \widetilde S}\prod_{i, j}\left | 2\sin \left(\frac{\tilde s-\tilde t}{r+k}+\frac{i-j}{\delta}\right)\pi \right |.$$ Here, we wrote $$s=\tilde s+\frac{r+k}{\delta}i\in \widetilde S_{i},\,\, t=\tilde t+\frac{r+k}{\delta}j\in \widetilde S_{j},$$ where $\tilde s, \tilde t\in \widetilde S$, $0\leq i,j\leq \delta-1$. To evaluate the product, we repeatedly use the identity \begin{equation}\label{sins}\prod_{i=0}^{\delta-1}\left|2 \sin \left(x+\frac{i\pi}{\delta}\right)\right|=\left|2\sin (\delta x)\right|.\end{equation} In particular, $$\prod_{0\leq i\neq j\leq \delta-1}\left|2\sin \frac{i-j}{\delta}\pi\right|=\delta^{\delta}.$$ We obtain \begin{equation}\label{deltah}\prod_{s\neq t\in S}\left | 2\sin \frac{s-t}{r+k}\pi \right |=\delta^{r}\prod_{\tilde s \neq \tilde t\in \widetilde S}\left |2 \sin \left(\frac{(\tilde s-\tilde t)\delta}{r+k}\pi\right)\right|^{\delta}.\end{equation} The prefactor $\delta^{r}$ accounts for the $\frac{r}{\delta}$ cases when $\tilde s=\tilde t$. 

From \eqref{sumbfn} we see that for any set $S$ we have \begin{equation} \label{xid}\xi_{d}(S)=\frac{1}{d^{2g}}\sum_{\delta|d\text { and } \delta|\delta_{S}}\bfn(\delta).\end{equation}
Equations \eqref{deltah} and \eqref{xid} allow us to match \eqref{comba} and the formula of Corollary \ref{combf} term by term, completing the proof. \qed

\end {document}